\documentclass{amsart}
\usepackage{graphicx}
\usepackage[subpreambles=true]{standalone}
\usepackage{amsmath,amssymb}
\usepackage{amsthm}
\usepackage{tikz}
\usetikzlibrary{calc}
\usepackage[active]{srcltx}
\usepackage{color}
\usepackage{graphicx}
\usepackage{xcolor}
\usepackage{enumitem}
\usepackage[margin=25mm]{geometry}
\usepackage{subcaption}
\renewcommand{\L}{\mathcal{L}}
\newcommand{\D}{\mathcal{D}}
\newcommand{\C}{\mathcal{C}}
\newcommand{\T}{\mathcal{T}}
\theoremstyle{plain}
\newtheorem{thm}{Theorem}[section]
\newtheorem{lem}{Lemma}[section]

\newtheorem{que}{Question}
\theoremstyle{definition}
\newtheorem{dfn}{Definition}[section]

\newtheorem*{theorem}{Main Theorem}
\newtheorem*{ack}{Acknowledgement}

\begin{document}

\setlength{\baselineskip}{12pt}

\title{Estimating the Kirby-Thompson invariants of surface-links in the Yoshikawa table}
\author{Minami Taniguchi}

\keywords{Knotted surface, Bridge trisection, Kirby-Thompson invariant}

\subjclass[2020]{Primary 57K10; Secondary 57K20}

\begin{abstract}
In \cite{tri_plane_diagrams}, minimal tri-plane diagrams of surface-links listed in the Yoshikawa table are computed using ch-diagrams. In this paper, we obtain upper bounds for the $\L$- and $\L^{\star}$-invariants of several surface-links in the Yoshikawa table by using their tri-plane diagrams.
\end{abstract}

\maketitle

\section{Introduction}

In \cite{bridge_trisection_1}, Meier--Zupan introduced a \textit{bridge trisection} of a surface-link $F$ in $S^{4}$, which is a decomposition of $(S^{4},F)$ into three trivial disk systems. In this way, we obtain a new diagram of $F$ called a \textit{tri-plane diagram}, which is a tuple of planar diagrams of three trivial tangles representing the boundaries of the trivial disk systems of the bridge trisection. In \cite{tri_plane_diagrams}, minimal tri-plane diagrams of all surface-links listed in the Yoshikawa table are computed, and their bridge numbers are specified.

In \cite{Kirby_Thompson_1}, Blair--Campisi--Taylor--Tomova defined the \textit{$\L$-invariant} $\L(F)$ of a surface-link $F\subset{S^{4}}$ using minimal bridge trisections of $F$ and the pants complex of the bridge surface of the bridge trisections, and they showed that any bridge trisection of a surface-link $F$ is standard if $\L(F)=0$. Using the dual curve complex of the bridge surface instead of the pants complex, Aranda--Pongtanapaisan--Zhang introduced the \textit{$\L^{\star}$-invariant} $\L^{\star}(F)$, and they showed that any bridge trisection of a surface-link $F$ is standard if $\L^{\star}(F)\leq2$. In other words, if $\L(F)=0$ or $\L^{\star}(F)\leq2$, a surface-link $F$ is smoothly trivial. We refer to both the $\L$- and $\L^{\star}$-invariants as the \textit{Kirby-Thompson invariants}, and these invariants measure the complexity of surface-links. On the other hand, there are a few examples of surface-links determined their values of the Kirby-Thompson invariants since it is not even straightforward to calculate upper bounds of these invariants because their definitions are so complicated. 

In this article, we compute upper bounds for the Kirby-Thompson invariants of several surface-links listed in the Yoshikawa table. 

\begin{theorem}[Theorem 3.1]
The Kirby-Thompson invariants of several surface-links in the Yoshikawa table are determined or estimated as follows.
\vspace{2mm}
\begin{center}
\begin{tabular}{|c||c|c|c|c|c|c|c|c|c|} \hline
\text{label} & $6_{1}^{0,1}$ & $7_{1}^{0,-2}$ & $8_{1}^{-1,-1}$ & $9_{1}$ & $9_{1}^{1,-2}$ & $10_{1}^{1}$ & $10_{3}$ & $10_{1}^{0,0,1}$ & $10_{1}^{-2,-2}$ \rule[-5pt]{0pt}{17pt} \\ \hline\hline
$\L$ & $15$ & $15\sim16$ & $15\sim18$ & $15\sim16$ & $\sim25$ & $\sim24$ & $15\sim34$ & $\sim27$ & $\sim26$ \rule[-5pt]{0pt}{17pt} \\ \hline
$\L^{\star}$ & $12$ & $12\sim13$ & $12$ & $12\sim13$ & $\sim19$ & $\sim18$ & $12\sim34$ & $\sim21$ & $\sim20$ \rule[-5pt]{0pt}{17pt} \\ \hline
\end{tabular}
\end{center}
\vspace{2mm}
\end{theorem}

In \cite{Kirby_Thompson_2} and \cite{Kirby_Thompson_3}, the values of the invariants of surface-links $8_{1}$, $8_{1}^{1,1}$, $10_{1}$, and $10_{2}$ are determined.

\begin{que}
Compute upper bounds for the Kirby-Thompson invariants of the other surface-links $9_{1}^{0,1}$, $10_{1}^{0,1}$, $10_{2}^{0,1}$, $10_{1}^{1,1}$, $10_{1}^{0,-2}$, $10_{2}^{0,-2}$, and $10_{1}^{-1,-1}$.
\end{que}

\begin{ack}

The author would like to express sincere gratitude to his supervisor, Hisaaki Endo, for his support and encouragement throughout this work.

\end{ack}

\section{Preliminaries}

In this section, we recall some concepts to define and calculate the Kirby-Thompson invariants of surface-links. First of all, we recall bridge trisections of surface-links in $S^{4}$. Let $F$ be a surface-link in $S^{4}$, and $X_{1}\cup X_{2}\cup X_{3}$ be a genus zero trisection of $S^{4}$. For each $i,j\in\{1,2,3\}$, let $(X_{i},\D_{i}):=(X_{i},X_{i}\cap F)$, $(B_{ij},\alpha_{ij}):=(X_{i},\D_{i})\cap(X_{j},\D_{j})=(\partial X_{i},\partial\D_{i})\cap(\partial X_{j},\partial\D_{j})$, and $(\Sigma,\mathrm{p}):=(X_{1},\D_{1})\cap(X_{2},\D_{2})\cap(X_{3},\D_{3})$. 

\begin{dfn}
The decomposition $(S^{4},F)=(X_{1},\D_{1})\cup(X_{2},\D_{2})\cup(X_{3},\D_{3})$ is said to be a \textit{$(b;c_{1},c_{2},c_{3})$-bridge trisection} if it satisfies the following properties.
\begin{enumerate}[label=(\roman*)]
\item $(X_{i},\D_{i})$ is a trivial $c_{i}$-disk system for each $i\in\{1,2,3\}$.
\item $(B_{ij},\alpha_{ij})$ is a trivial $b$-tangle for each $i\not=j\in\{1,2,3\}$.
\item $(\Sigma,\mathrm{p})$ is a $2b$-punctured 2-sphere, where $\mathrm{p}=\alpha_{ij}\cap\Sigma=\alpha_{ij}\cap\alpha_{ki}$.
\end{enumerate}
\end{dfn}

\begin{figure}[t]
\centering
\input{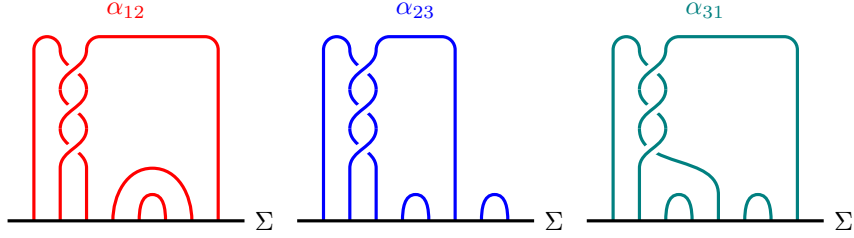} 
\caption{A tri-plane diagram of spun trefoil.}
\label{fig:tri-plane diagram of spun trefoil}
\end{figure}

Let $(S^{4},F)=(X_{1},\D_{1})\cup(X_{2},\D_{2})\cup(X_{3},\D_{3})$ be a $(b;c_{1},c_{2},c_{3})$-bridge trisection. From the above definition, each $(B_{ij},\alpha_{ij})$ is a trivial $b$-tangle, and the boundary sum $(B_{ij},\alpha_{ij})\cup_{\partial}\overline{(B_{ki},\alpha_{ki})}$ is the unlink $(\partial X_{i},\partial\D_{i})$. We call the tuple of trivial $b$-tangles as a \textit{tri-plane diagram} of the bridge trisection. For instance, a tri-plane diagram of the spun trefoil is shown in Figure \ref{fig:tri-plane diagram of spun trefoil}. Refer to \cite{Kirby_Thompson_2} for more detailed information on how to depict the diagram, and to \cite{banded_unlink_diagram} for more detailed information on banded unlink diagrams. 

%tri-plane diagram

Secondly, we recall the definitions of the pants complex and the dual curve complex of a punctured 2-sphere. Let $b$ be an integer with $b\geq2$. Let $\Sigma=\Sigma_{2b,0}$ be a $2b$-punctured 2-sphere, and $\C=\{c_{1},\cdots,c_{2b-3}\}$ be a collection of disjoint essential simple closed curves on $\Sigma$.

%A-move

\begin{figure}[b]
\centering
\input{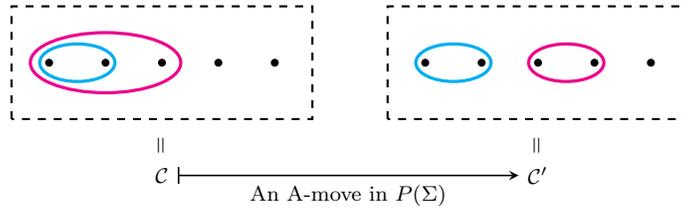} 
\caption{An example of pants decompositions $\C$ and $\C^{\prime}$ of five punctured 2-sphere $\Sigma_{5,\,0}$ whose distance in $P(\Sigma)$ is equal to one.}
\end{figure}

\begin{dfn}
\phantom{a}
\begin{enumerate}[label=(\roman*)]
\item $\C$ is said to be a \textit{pants decomposition} of $\Sigma$ if we obtain $2b-2$ pairs of pants by cutting $\Sigma$ along every essential curve $c\in\C$.
\item Let $\C=\{c_{1},\cdots,c_{2b-3}\}$ and $\C^{\prime}=\{c_{1}^{\prime},\cdots,c_{2b-3}^{\prime}\}$ be pants decompositions of $\Sigma$. $\C^{\prime}$ is obtained from $\C$ by an \textit{A-move} if they satisfy the following properties.
\begin{enumerate}[label=(\Roman*)]
\item For each $i\in\{1,\cdots,2b-2\}$, $c_{i}=c_{i}^{\prime}$.
\item $|c_{2b-3}\cap c_{2b-3}^{\prime}|=2$.
\end{enumerate}
\vspace{1mm}
Similarly, $\C^{\prime}$ is obtained from $\C$ by an \textit{$\text{A}^{\star}$-move} if they satisfy the following properties.
\begin{enumerate}[label=(\Roman*)]
\item For each $i\in\{1,\cdots,2b-2\}$, $c_{i}=c_{i}^{\prime}$.
\item $|c_{2b-3}\cap c_{2b-3}^{\prime}|\geq2$.
\end{enumerate}
\vspace{1mm}
\item A 1-complex $P(\Sigma)$ is said to be the \textit{pants complex} of $\Sigma$ if it satisfies the following properties.
\begin{enumerate}[label=(\Roman*)]
\item A vertex of $P(\Sigma)$ is a pants decomposition $\C$ of $\Sigma$.
\item Two vertex $\C,\C^{\prime}\in P(\Sigma)$ is connected by an \textit{edge} in $P(\Sigma)$ if $\C^{\prime}$ is obtained from $\C$ by an A-move.
\end{enumerate}
Similarly, a 1-complex $C^{\star}(\Sigma)$ is said to be the \textit{dual curve complex} if it satisfies the following properties.
\begin{enumerate}[label=(\Roman*)]
\item A vertex of $C^{\star}(\Sigma)$ is a pants decomposition $\C$ of $\Sigma$.
\item Two vertices $\C,\C^{\prime}\in{C^{\star}(\Sigma)}$ is connected by an \textit{edge} in $C^{\star}(\Sigma)$ if $\C^{\prime}$ is obtained from $\C$ by an $\text{A}^{\star}$-move.
\end{enumerate}
\end{enumerate}
\end{dfn}

Let $(B,\alpha)$ be a trivial $b$-tangle and $(\Sigma,\mathrm{p})=(\partial B,\partial\alpha)$ the $2b$-punctured 2-sphere. Let $c$ be an essential simple closed curve on $\Sigma$. 

\begin{dfn}
\phantom{a}
\begin{enumerate}[label=(\roman*)]
\item The curve $c$ is said to be \textit{compressing} if there exists a 2-disk $D\subset{B}$ such that $\partial D=c$ and $|D\cap\alpha|=0$. We refer to such a 2-disk $D$ as a \textit{compressing disk} for the tangle $(B,\,\alpha)$.
\item The curve $c$ is said to be \textit{cut} if there exists a 2-disk $D\subset{B}$ such that $\partial D=c$ and $|D\cap\alpha|=1$. We refer to such a 2-disk $D$ as a \textit{cut disk} for the tangle $(B,\alpha)$.
\end{enumerate}
\end{dfn}

%c-reducing

\begin{figure}[t]
\centering
\begin{subfigure}[b]{0.4\textwidth}
\centering
\input{c_reducing.tex}
\caption{An example of a compressing curve $c$ and a cut curve $c^{\prime}$ for a trivial 3-tangle $\alpha$.}
\end{subfigure}
\hspace{1mm}
\begin{subfigure}[b]{0.4\textwidth}
\centering
\input{c_reducing_2.tex}
\caption{An example of a reducing curve $c$ and a c-reducing curve $c^{\prime}$ for a 3-bridge split unlink $L$.}
\end{subfigure}
\end{figure}

Let $(S^{3},L)$ be a $b$-bridge split unlink $(B_{+},\alpha_{+})\cup_{\partial}\overline{(B_{-},\alpha_{-})}$, where each $(B_{\pm},\alpha_{\pm})$ is a trivial $b$-tangle, and $(\Sigma,\mathrm{p})=(\partial B_{\pm},\partial\alpha_{\pm})$ the $2b$-punctured 2-sphere. Let $c$ be an essential simple closed curve on $\Sigma$.

\begin{dfn}
\phantom{a}
\begin{enumerate}[label=(\roman*)]
\item The curve $c$ is said to be \textit{reducing} if there exists a 2-sphere $Q\subset{S^{3}}$ such that $Q\cap\Sigma=c$ and $|Q\cap L|=0$. We refer to such a 2-sphere $Q$ as a \textit{reducing sphere} for the unlink $(S^{3},L)$. We note that each 2-disk $Q\cap B_{\pm}\subset{B_{\pm}}$ is also compressing for the tangle $(B_{\pm},\alpha_{\pm})$, respectively.
\item The curve $c$ is said to be \textit{c-reducing} if there exists a 2-sphere $Q\subset{S^{3}}$ such that $Q\cap\Sigma=c$ and $|Q\cap\alpha_{+}|=|Q\cap\alpha_{+}|=1$. We refer to such a 2-sphere $Q$ as a \textit{c-reducing sphere} for the unlink $(S^{3},L)$. We note that each 2-disk $Q\cap B_{\pm}\subset{B_{\pm}}$ is also cut for the tangle $(B_{\pm},\alpha_{\pm})$, respectively.
\end{enumerate}
\end{dfn}

With these preliminaries, we define the Kirby-Thompson invariants of a surface-link. Let $F\subset{S^{4}}$ be a surface-link and $(S^{4},F)=(X_{1},\D_{1})\cup(X_{2},\D_{2})\cup(X_{3},\D_{3})$ a minimal $(b;c_{1},c_{2},c_{3})$-bridge trisection of $F$. Let $\T$ denote the bridge trisection. Let $(\Sigma,\mathrm{p})=(X_{1},\D_{1})\cap(X_{2},\D_{2})\cap(X_{3},\D_{3})$ be the bridge surface of the bridge trisection, and we regard the bridge surface as a $2b$-punctured 2-sphere. Let $D_{c}(\alpha_{ij})$ denote the set of vertices of $P(\Sigma)$ whose curves are either compressing or cut curves for the trivial $b$-tangle $(B_{ij},\alpha_{ij})$.

\begin{figure}[htbp]
\centering
\input{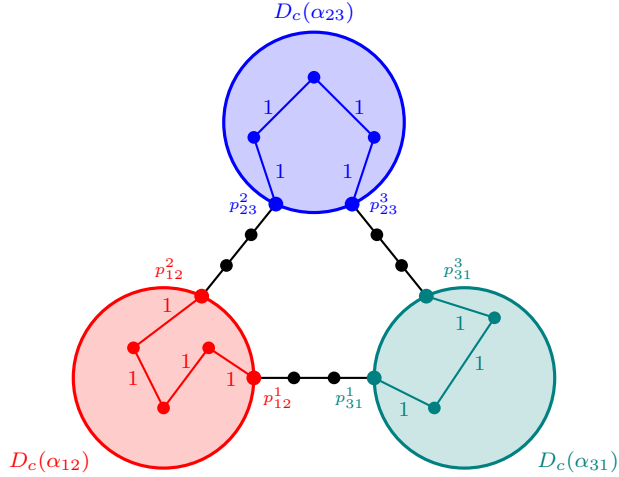}
\caption{An image of the Kirby-Thompson invariants, where each $(p_{ij}^{i},p_{ki}^{i})$ is an efficient defining pair for the unlink $\partial\D_{i}$.}
\end{figure}

\begin{dfn}
\phantom{a}
\begin{enumerate}[label=(\roman*)]
\item A pair of efficient pants decompositions $(p_{ij}^{i},p_{ki}^{i})\in{D_{c}(\alpha_{ij})\times D_{c}(\alpha_{ki})}$ is said to be an \textit{efficient defining pair} for the unlink $\partial\D_{i}$ if the distance between $p_{ij}^{i}$ and $p_{ki}^{i}$ in $P(\Sigma)$ (or in $C^{\star}(\Sigma)$) is equal to $b-c_{i}$.
\item We define the \textit{$\L$-invariant} $\L(\T)$ of the bridge trisection $\T$ as \[
\L(\T):=\min\left\{\sum_{i,\,j}d^{P}(p_{ij}^{i},p_{ij}^{j})\,\middle|\,
\begin{tabular}{l}
\text{$(p_{ij}^{i},p_{ki}^{i})\in{D_{c}(\alpha_{ij})\times D_{c}(\alpha_{ki})}$ is} \\
\text{an efficient defining pair for $\partial\D_{i}$.}
\end{tabular}
\right\},
\]
where $d^{P}\colon P(\Sigma)\times P(\Sigma)\longrightarrow \mathbb{Z}_{\geq0}$ is the metric function of $P(\Sigma)$.
\item We define the \textit{$\L$-invariant} $\L(F)$ of the surface-link $F$ as follows. \[
\L(F):=\min\left\{\L(\T)\,\middle|\,\text{$\T$ is a $b$-bridge trisection of $F$.}\right\}.
\]
\end{enumerate}
Using the dual curve complex $C^{\star}(\Sigma)$ instead of the pants complex $P(\Sigma)$, we can define another invariant $\L^{\star}$ as follows.
\begin{enumerate}[label=(\roman*),start=4]
\item We define the \textit{$\L^{\star}$-invariant} $\L^{\star}(\T)$ of the bridge trisection $\T$ as \[
\L^{\star}(\T):=\min\left\{\sum_{i,\,j}d^{\star}(p_{ij}^{i},p_{ij}^{j})\,\middle|\,
\begin{tabular}{l}
\text{$(p_{ij}^{i},p_{ki}^{i})\in{D_{c}(\alpha_{ij})\times D_{c}(\alpha_{ki})}$ is} \\
 \text{an efficient defining pair for $\partial\D_{i}$.}
 \end{tabular}
 \right\},
\]
where $d^{\star}\colon C^{\star}(\Sigma)\times C^{\star}(\Sigma)\longrightarrow\mathbb{Z}_{\geq0}$ is the metric function of $C^{\star}(\Sigma)$.
\item We define the \textit{$\L^{\star}$-invariant} $\L^{\star}(F)$ of the surface-link $F$ as follows. \[
\L^{\star}(F):=\min\left\{\L^{\star}(\T)\,\middle|\,\text{$\T$ is a $b$-bridge trisection of $F$.}\right\}.
\]
\end{enumerate}
We refer to both the $\L$-invariant and the $\L^{\star}$-invariant as the \textit{Kirby-Thompson invariants}.
\end{dfn}

%図4

The following properties of the Kirby-Thompson invariants are fundamental.

\begin{thm}[\cite{Kirby_Thompson_1}, \cite{Kirby_Thompson_3}]
Let $F$ be a surface-link in $S^{4}$, then the following properties hold.
\begin{enumerate}[label=(\roman*)]
\item $\L^{\star}(F),\,\L(F)\in\mathbb{Z}_{\geq0},\,\L^{\star}(F)\leq\L(F)$.
\item If $\L(F)=0$ holds, the surface-link $F$ is smoothly isotopic to the distant sum of the connected sum of finitely many standard surfaces.
\item If $\L^{\star}(F)\leq2$ holds, the surface-link $F$ is smoothly isotopic to the distant sum of the connected sum of finitely many standard surfaces.
\end{enumerate}
\end{thm}

In other words, the Kirby-Thompson invariants can specify whether a surface-link is smoothly unlinked or not. Furthermore, in \cite{Kirby_Thompson_3}, some lower bounds for these invariants are provided for irreducible and unstabilized bridge trisections.

\begin{lem}[\cite{Kirby_Thompson_3}]
Let $\T$ be an irreducible and unstabilized $(b;c)$-bridge trisection for a surface-link $F$. Then,
\begin{enumerate}[label=(\roman*)]
\item $\L(\T)\geq\L^{\star}(\T)\geq3(b+c-3)$.
\item $(b,c)=(4,2)\,\Longrightarrow\,\L^{\star}(\T)\geq12$.
\item $c\geq2\,\Longrightarrow\,\L(\T)\geq3(b+c-2)$.
\item $c=2\,\Longrightarrow\,\L(\T)\geq3(b+1)$.
\end{enumerate}
\label{lem:lower}
\end{lem}

\section{Main result}

In this section, we prove the following main theorem. We note that, in all of the following figures, all compressing or cut curves for tangles $\alpha_{12}$, $\alpha_{23}$, and $\alpha_{31}$ are depicted in red, blue, and green, respectively, and all common curves in $p_{ij}^{i}\cap p_{ki}^{k}$ are depicted in orange for each $\{i,j,k\}=\{1,2,3\}$.

\begin{thm}
The Kirby-Thompson invariants of several surface-links in the Yoshikawa table are determined or estimated as follows.
\vspace{2mm}
\begin{center}
\begin{tabular}{|c||c|c|c|c|c|c|c|c|c|} \hline
\text{label} & $6_{1}^{0,1}$ & $7_{1}^{0,-2}$ & $8_{1}^{-1,-1}$ & $9_{1}$ & $9_{1}^{1,-2}$ & $10_{1}^{1}$ & $10_{3}$ & $10_{1}^{0,0,1}$ & $10_{1}^{-2,-2}$ \rule[-5pt]{0pt}{17pt} \\ \hline\hline
$\L$ & $15$ & $15\sim16$ & $15\sim18$ & $15\sim16$ & $\sim25$ & $\sim24$ & $15\sim34$ & $\sim27$ & $\sim26$ \rule[-5pt]{0pt}{17pt} \\ \hline
$\L^{\star}$ & $12$ & $12\sim13$ & $12$ & $12\sim13$ & $\sim19$ & $\sim18$ & $12\sim34$ & $\sim21$ & $\sim20$ \rule[-5pt]{0pt}{17pt} \\ \hline
\end{tabular}
\end{center}
\vspace{2mm}
\label{thm:main}
\end{thm}

\begin{proof}
From lemma \ref{lem:6_1^01}, lemma \ref{lem:7_1^0-2}, lemma \ref{lem:8_1^-1-1}, lemma \ref{lem:9_1}, lemma \ref{lem:9_1^1-2}, lemma \ref{lem:10_1^1}, lemma \ref{lem:10_3}, lemma \ref{lem:10_1^001}, and lemma \ref{lem:10_1^-2-2}, we obtain upper bounds for the Kirby-Thompson invariants of the surface-links listed in the above table. On the other hand, since any minimal $(4;2)$-bridge trisection is irreducible, from lemma \ref{lem:lower}, we obtain lower bounds for the Kirby-Thompson invariants of surface-links whose bridge number are four.
\end{proof}

In \cite{Kirby_Thompson_3}, any $(6;2)$-minimal bridge trisection of $T^{2}$-spin of nontrivial 2-bridge link is irreducible. If any minimal bridge trisection with bridge number less than or equal to six is irreducible, we obtain stricter evaluations for the invariants of the surface-links in theorem \ref{thm:main}.

\begin{que}
Is any minimal bridge trisection with bridge number less than or equal to six irreducible?
\end{que}

\begin{lem}
The following inequalities hold. \[
\L(6_{1}^{0,1})\leq15,\,\L^{\star}(6_{1}^{0,1})\leq12
\]
\label{lem:6_1^01}
\end{lem}

\begin{figure}[htbp]
\centering
\input{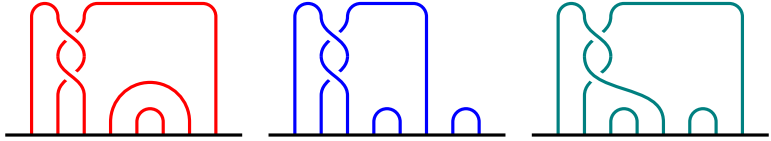}
\caption{A tri-plane diagram of $6_{1}^{0,1}$.}
\label{fig:6_1^01}
\end{figure}

\begin{proof}
From \cite{tri_plane_diagrams}, a tri-plane diagram of $6_{1}^{0,1}$ is as in Figure \ref{fig:6_1^01}. From the tri-plane diagram in Figure \ref{fig:6_1^01}, we can find efficient defining pairs for the bridge trisection of $6_{1}^{0,1}$ shown in Figure \ref{fig:L_6_1^01}. We can see that \[
d^{P}(p_{ij}^{i},p_{ij}^{j})\leq5
\]
for each $i,j\in\{1,2,3\}$, hence $\L(6_{1}^{0,1})\leq5+5+5=15$ holds. In the same way, from Figure \ref{fig:L*_6_1^01}, we see that \[
 d^{\star}(p_{ij}^{i},p_{ij}^{j})\leq4
 \]
 for each $i,j\in\{1,2,3\}$, hence $\L^{\star}(6_{1}^{0,1})\leq4+4+4=12$ holds.
\end{proof}

\begin{lem}
The following inequalities hold. \[
\L(7_{1}^{0,-2})\leq16,\,\L^{\star}(7_{1}^{0,-2})\leq13
\]
\label{lem:7_1^0-2}
\end{lem}

\begin{figure}[htbp]
\centering
\input{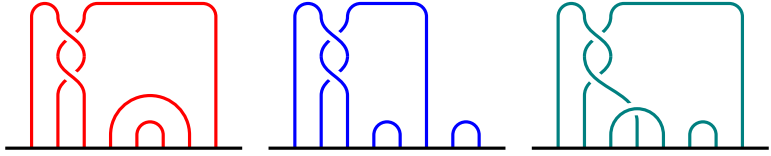}
\caption{A tri-plane diagram of $7_{1}^{0,-2}$.}
\label{fig:7_1^0-2}
\end{figure}

\begin{proof}
From \cite{tri_plane_diagrams}, a tri-plane diagram of $7_{1}^{0,-2}$ is as in Figure \ref{fig:7_1^0-2}. From the tri-plane diagram in Figure \ref{fig:7_1^0-2}, we can find efficient defining pairs for the bridge trisection of $7_{1}^{0,-2}$ shown in Figure \ref{fig:L_7_1^0-2}. We can see that \[
d^{P}(p_{12}^{1},p_{12}^{2})\leq4,\,d^{P}(p_{23}^{2},p_{23}^{3})\leq4,\,d^{P}(p_{31}^{3},p_{31}^{1})\leq5
\]
for each $i,j\in\{1,2,3\}$, hence $\L(7_{1}^{0,-2})\leq4+4+5=13$ holds. In the same way, from Figure \ref{fig:L*_7_1^0-2}, we see that \[
 d^{\star}(p_{12}^{1},p_{12}^{2})\leq4,\,d^{\star}(p_{23}^{2},p_{23}^{3})\leq4,\,d^{\star}(p_{31}^{3},p_{31}^{1})\leq5,
 \]
 hence $\L^{\star}(7_{1}^{0,-2})\leq4+4+5=13$ holds.
\end{proof}

\begin{que}
Does there exist some bridge trisection of $7_{1}^{0,-2}$ together with efficient defining pairs such that $\L(7_{1}^{0,-2})\leq15,\,\L^{\star}(7_{1}^{0,-2})\leq12$?
\end{que}

\begin{lem}
The following inequalities hold. \[
\L(8_{1}^{-1,-1})\leq18,\,\L^{\star}(8_{1}^{-1,-1})\leq12
\]
\label{lem:8_1^-1-1}
\end{lem}

\begin{figure}[htbp]
\centering
\input{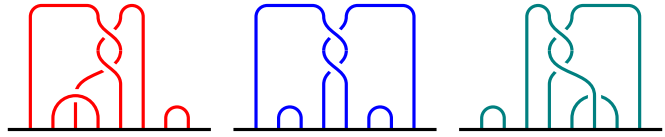}
\caption{A tri-plane diagram of $8_{1}^{-1,-1}$.}
\label{fig:8_1^-1-1}
\end{figure}

\begin{proof}
From \cite{tri_plane_diagrams}, a tri-plane diagram of $8_{1}^{-1,-1}$ is as in Figure \ref{fig:8_1^-1-1}. From the tri-plane diagram in Figure \ref{fig:8_1^-1-1}, we can find efficient defining pairs for the bridge trisection of $8_{1}^{-1,-1}$ shown in Figure \ref{fig:L_8_1^-1-1}. We see that \[
d^{P}(p_{ij}^{i},p_{ij}^{j})\leq6
\]
 for each $i,j\in\{1,2,3\}$, hence $\L(8_{1}^{-1,-1})\leq6+6+6=18$ holds. In the same way, from Figure \ref{fig:L*_8_1^-1-1}, we see that \[
 d^{\star}(p_{ij}^{i},p_{ij}^{j})\leq4
 \]
 for each $i,j\in\{1,2,3\}$, hence $\L^{\star}(8_{1}^{-1,-1})\leq4+4+4=12$ holds.
\end{proof}

\begin{lem}
The following inequalities hold. \[
\L(9_{1})\leq16,\,\L^{\star}(9_{1})\leq13
\]
\label{lem:9_1}
\end{lem}

\begin{proof}
From \cite{tri_plane_diagrams}, a tri-plane diagram of $9_{1}$ is as in Figure \ref{fig:9_1}. From the tri-plane diagram in Figure \ref{fig:9_1}, we can find efficient defining pairs for the bridge trisection of $9_{1}$ shown in Figure \ref{fig:L_9_1}. We see that \[
d^{P}(p_{12}^{1},p_{12}^{2})\leq5,\,d^{P}(p_{23}^{2},p_{23}^{3})\leq5,\,d^{P}(p_{31}^{3},p_{31}^{1})\leq6,
\]
 hence $\L(9_{1})\leq5+5+6=16$ holds. In the same way, from Figure \ref{fig:L*_9_1}, we see that \[
 d^{\star}(p_{12}^{1},p_{12}^{2})\leq4,\,d^{\star}(p_{23}^{2},p_{23}^{3})\leq4,\,d^{\star}(p_{31}^{3},p_{31}^{1})\leq5,
 \]
 hence $\L^{\star}(9_{1})\leq4+4+5=13$ holds. 
\end{proof}

\begin{que}
Does there exist some bridge trisection of $9_{1}$ together with efficient defining pairs such that $\L(9_{1})\leq15,\,\L^{\star}(9_{1})\leq12$?
\end{que}

\begin{lem}
The following inequations hold. \[
\L(9_{1}^{1,-2})\leq25,\,\L^{\star}(9_{1}^{1,-2})\leq19
\]
\label{lem:9_1^1-2}
\end{lem}

\begin{proof}
From \cite{tri_plane_diagrams}, a tri-plane diagram of $9_{1}^{1,-2}$ is as in Figure \ref{fig:9_1^1-2}. From the tri-plane diagram in Figure \ref{fig:9_1^1-2}, we can find efficient defining pairs for the bridge trisection of $9_{1}^{1,-2}$ shown in Figure \ref{fig:L_9_1^1-2}. We see that \[
d^{P}(p_{12}^{1},p_{12}^{2})\leq9,\,d^{P}(p_{23}^{2},p_{23}^{3})\leq8,\,d^{P}(p_{31}^{3},p_{31}^{1})\leq8, 
\]
hence $\L(9_{1}^{1,-2})\leq9+8+8=25$ holds. In the same way, from Figure \ref{fig:L*_9_1^1-2}, we see that \[
d^{\star}(p_{12}^{1},p_{12}^{2})\leq7,\,d^{\star}(p_{23}^{2},p_{23}^{3})\leq6,\,d^{\star}(p_{31}^{3},p_{31}^{1})\leq6, 
\]
hence $\L^{\star}(9_{1}^{1,-2})\leq7+6+6=19$ holds.
\end{proof}

\begin{que}
Does there exist some bridge trisection of $9_{1}^{1,-2}$ together with efficient defining pairs such that $\L(9_{1}^{1,-2})\leq24,\,\L^{\star}(9_{1}^{1,-2})\leq18$?
\end{que}

\begin{lem}
The following inequations hold. \[
\L(10_{1}^{1})\leq24,\,\L^{\star}(10_{1}^{1})\leq18
\]
\label{lem:10_1^1}
\end{lem}

\begin{proof}
From \cite{tri_plane_diagrams}, a tri-plane diagram of $10_{1}^{1}$ is as in Figure \ref{fig:10_1^1}. From the tri-plane diagram in Figure \ref{fig:10_1^1}, we can find efficient defining pairs for the bridge trisection of $10_{1}^{1}$ shown in Figure \ref{fig:L_10_1^1}. We see that \[
d^{P}(p_{ij}^{i},p_{ij}^{j})\leq8
\]
for each $i,j\in\{1,2,3\}$, hence $\L(10_{1}^{1})\leq8+8+8=24$ holds. In the same way, from Figure \ref{fig:L*_10_1^1}, we see that \[
d^{\star}(p_{ij}^{i},p_{ij}^{j})\leq6
\]
for each $i,j\in\{1,2,3\}$, hence $\L^{\star}(10_{1}^{1})\leq6+6+6=18$ holds.

\end{proof}

\begin{figure}[htbp]
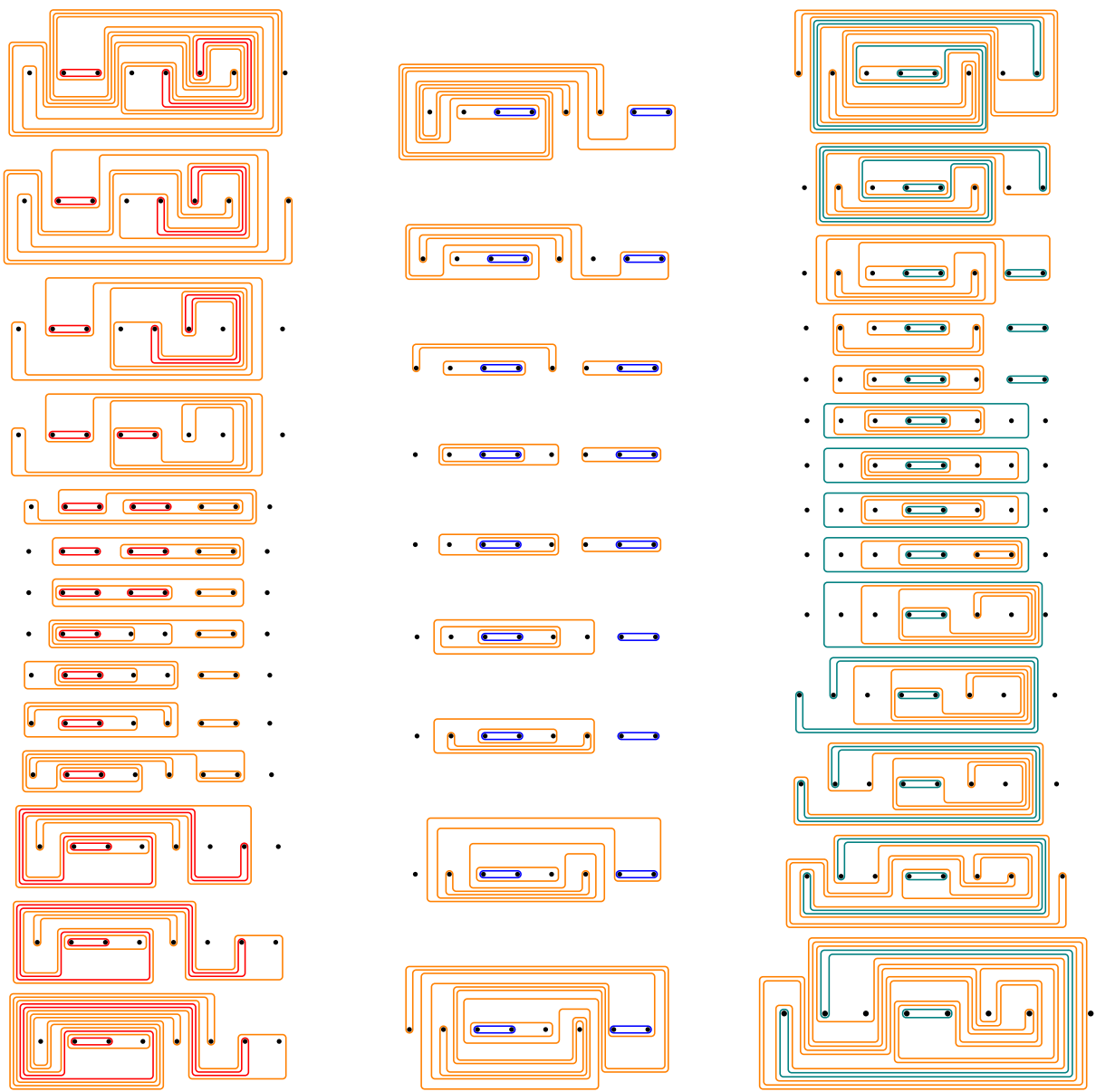

\centering
\begin{subfigure}[b]{0.3\textwidth}
\centering
\input{L_10_3_1.tex}
\end{subfigure}
\hspace{5mm}
\begin{subfigure}[b]{0.3\textwidth}
\centering
\input{L_10_3_2.tex}
\end{subfigure}
\hspace{5mm}
\begin{subfigure}[b]{0.3\textwidth}
\centering
\input{L_10_3_3.tex}
\end{subfigure}
\caption{An upper bound for the $\L$-invariant and the $\L^{\star}$-invariant of $10_{3}$.}
\label{fig:L_10_3}
\end{figure}

The efficient defining pairs for the minimal tri-plane diagram of $10_{3}$, constructed in Figure \ref{fig:L_10_3}, are so complicated that we cannot find any shorter paths in $C^{\star}(\Sigma)$. 

\begin{lem}
The following inequations hold. \[
\L^{\star}(10_{3})\leq\L(10_{3})\leq34
\]
\label{lem:10_3}
\end{lem}

\begin{figure}[htbp]
\centering
\input{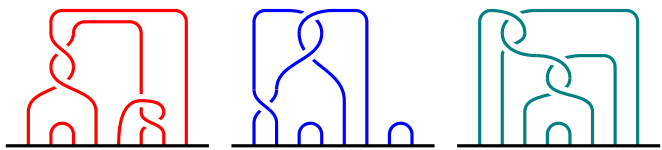}
\caption{A tri-plane diagram of $10_{3}$.}
\label{fig:10_3}
\end{figure}

\begin{proof}
From \cite{tri_plane_diagrams}, a tri-plane diagram of $10_{3}$ is as in Figure \ref{fig:10_3}. From the tri-plane diagram in Figure \ref{fig:10_3}, we can find efficient defining pairs for the bridge trisection of $10_{3}$ shown in Figure \ref{fig:L_10_3}. We see that \[
d^{P}(p_{12}^{1},p_{12}^{2})\leq13,\,d^{P}(p_{23}^{2},p_{23}^{3})\leq8,\,d^{P}(p_{31}^{3},p_{31}^{1})\leq13,
\]
hence $\L(10_{3})\leq13+8+13=34$ holds. 
\end{proof}

\begin{que}
Find a sequence of A-moves realizing minimal distance of $p_{ij}^{i}\rightarrow p_{ij}^{j}$ in $C^{\star}(\Sigma)$ for each $i\not=j\in\{1,2,3\}$. Furthermore, find alternative efficient defining pairs for the bridge trisection of $10_{3}$ that give stricter upper bounds of $\L(10_{3})$ and $\L^{\star}(10_{3})$.
\end{que}

\begin{lem}
The following inequations hold. \[
\L(10_{1}^{0,0,1})\leq27,\,\L^{\star}(10_{1}^{0,0,1})\leq21
\]
\label{lem:10_1^001}
\end{lem}

\begin{proof}
From \cite{tri_plane_diagrams}, a tri-plane diagram of $10_{1}^{0,0,1}$ is as in Figure \ref{fig:10_1^001}. From the tri-plane diagram in Figure \ref{fig:10_1^001}, we can find efficient defining pairs for the bridge trisection of $10_{1}^{0,0,1}$ shown in Figure \ref{fig:L_10_1^001}. We see that \[
d^{P}(p_{ij}^{i},p_{ij}^{j})\leq9
\]
for each $i,j\in\{1,2,3\}$, hence $\L(10_{1}^{0,0,1})\leq9+9+9=27$ holds. In the same way, from Figure \ref{fig:L*_10_1^001}, we see that \[
d^{\star}(p_{ij}^{i},p_{ij}^{j})\leq7
\]
for each $i,j\in\{1,2,3\}$, hence $\L^{\star}(10_{1}^{0,0,1})\leq7+7+7=21$ holds.
\end{proof}

\begin{lem}
The following inequalities hold. \[
\L(10_{1}^{-2,-2})\leq26,\,\L^{\star}(10_{1}^{-2,-2})\leq20
\]
\label{lem:10_1^-2-2}
\end{lem}

\begin{proof}
From \cite{tri_plane_diagrams}, a tri-plane diagram of $10_{1}^{-2,-2}$ is as in Figure \ref{fig:10_1^-2-2}. From the tri-plane diagram in Figure \ref{fig:10_1^-2-2}, we can find efficient defining pairs for the bridge trisection of $10_{1}^{-2,-2}$ shown in Figure \ref{fig:L_10_1^-2-2}. We see that \[
d^{P}(p_{12}^{1},p_{12}^{2})\leq10,\,d^{P}(p_{23}^{2},p_{23}^{3})\leq8,\,d^{P}(p_{31}^{3},p_{31}^{1})\leq8,
\]
hence $\L(10_{1}^{-2,-2})\leq10+8+8=26$ holds. In the same way, from Figure \ref{fig:L*_10_1^-2-2}, we see that \[
d^{\star}(p_{12}^{1},p_{12}^{2})\leq8,\,d^{\star}(p_{23}^{2},p_{23}^{3})\leq6,\,d^{\star}(p_{31}^{3},p_{31}^{1})\leq6,
\]
hence $\L^{\star}(10_{1}^{-2,-2})\leq8+6+6=20$ holds.
\end{proof}

\begin{que}
Does there exist some bridge trisection of $10_{1}^{-2,-2}$ together with efficient defining pairs such that $\L(10_{1}^{-2,-2})\leq24,\,\L^{\star}(10_{1}^{-2,-2})\leq18$?
\end{que}

These lemmas suggest that the Kirby-Thompson invariants of surface-links may take values that are not multiples of three, although all the values of the invariants specified so far are divisible by three. 

\begin{que}
 Do $\L(F)\equiv0,\,\L^{\star}(F)\equiv0\pmod3$ always hold for any surface-link $F$?
\end{que}

%========================================================

%=======================================================

\clearpage

\begin{figure}[htbp]
\centering
\begin{subfigure}[b]{0.3\textwidth}
\centering
\input{L_6_1_01.tex}
\end{subfigure}
\hspace{5mm}
\begin{subfigure}[b]{0.3\textwidth}
\centering
\input{L_6_1_01_2.tex}
\end{subfigure}
\hspace{5mm}
\begin{subfigure}[b]{0.3\textwidth}
\centering
\input{L_6_1_01_3.tex}
\end{subfigure}
\caption{An upper bound for the $\L$-invariant of $6_{1}^{0,1}$.}
\label{fig:L_6_1^01}
\end{figure}

\begin{figure}[htbp]
\centering
\begin{subfigure}[b]{0.3\textwidth}
\centering
\input{_L_6_1_01_1.tex}
\end{subfigure}
\hspace{5mm}
\begin{subfigure}[b]{0.3\textwidth}
\centering
\input{_L_6_1_01_2.tex}
\end{subfigure}
\hspace{5mm}
\begin{subfigure}[b]{0.3\textwidth}
\centering
\input{_L_6_1_01_3.tex}
\end{subfigure}
\caption{An upper bound for the $\L^{\star}$-invariant of $6_{1}^{0,1}$.}
\label{fig:L*_6_1^01}
\end{figure}

\clearpage

\begin{figure}[htbp]
\centering
\begin{subfigure}[b]{0.3\textwidth}
\centering
\input{L_7_1_0-2_1.tex}
\end{subfigure}
\hspace{5mm}
\begin{subfigure}[b]{0.3\textwidth}
\centering
\input{L_7_1_0-2_2.tex}
\end{subfigure}
\hspace{5mm}
\begin{subfigure}[b]{0.3\textwidth}
\centering
\input{L_7_1_0-2_3.tex}
\end{subfigure}
\caption{An upper bound for the $\L$-invariant of $7_{1}^{0,-2}$.}
\label{fig:L_7_1^0-2}
\end{figure}

\begin{figure}[htbp]
\centering
\begin{subfigure}[b]{0.3\textwidth}
\centering
\input{_L_7_1_0-2_1.tex}
\end{subfigure}
\hspace{5mm}
\begin{subfigure}[b]{0.3\textwidth}
\centering
\input{_L_7_1_0-2_2.tex}
\end{subfigure}
\hspace{5mm}
\begin{subfigure}[b]{0.3\textwidth}
\centering
\input{_L_7_1_0-2_3.tex}
\end{subfigure}
\caption{An upper bound for the $\L^{\star}$-invariant of $7_{1}^{0,-2}$.}
\label{fig:L*_7_1^0-2}
\end{figure}

\clearpage

\begin{figure}[htbp]
\centering
\begin{subfigure}[b]{0.3\textwidth}
\centering
\input{L_8_1_-1-1_1.tex}
\end{subfigure}
\hspace{5mm}
\begin{subfigure}[b]{0.3\textwidth}
\centering
\input{L_8_1_-1-1_2.tex}
\end{subfigure}
\hspace{5mm}
\begin{subfigure}[b]{0.3\textwidth}
\centering
\input{L_8_1_-1-1_3.tex}
\end{subfigure}
\caption{An upper bound for the $\L$-invariant of $8_{1}^{-1,-1}$.}
\label{fig:L_8_1^-1-1}
\end{figure}

\begin{figure}[htbp]
\centering
\begin{subfigure}[b]{0.3\textwidth}
\centering
\input{_L_8_1_-1-1_1.tex}
\end{subfigure}
\hspace{5mm}
\begin{subfigure}[b]{0.3\textwidth}
\centering
\input{_L_8_1_-1-1_2.tex}
\end{subfigure}
\hspace{5mm}
\begin{subfigure}[b]{0.3\textwidth}
\centering
\input{_L_8_1_-1-1_3.tex}
\end{subfigure}
\caption{An upper bound for the $\L^{\star}$-invariant of $8_{1}^{-1,-1}$.}
\label{fig:L*_8_1^-1-1}
\end{figure}

\clearpage

\begin{figure}[htbp]
\centering
\input{9_1.tex}
\caption{A tri-plane diagram of $9_{1}$.}
\label{fig:9_1}
\end{figure}

\begin{figure}[htbp]
\centering
\begin{subfigure}[b]{0.3\textwidth}
\centering
\input{L_9_1_1.tex}
\end{subfigure}
\hspace{5mm}
\begin{subfigure}[b]{0.3\textwidth}
\centering
\input{L_9_1_2}
\end{subfigure}
\hspace{5mm}
\begin{subfigure}[b]{0.3\textwidth}
\centering
\input{L_9_1_3.tex}
\end{subfigure}
\caption{An upper bound for the $\L$-invariant of $9_{1}$.}
\label{fig:L_9_1}
\end{figure}

\begin{figure}[htbp]
\centering
\begin{subfigure}[b]{0.3\textwidth}
\centering
\input{_L_9_1_1.tex}
\end{subfigure}
\hspace{5mm}
\begin{subfigure}[b]{0.3\textwidth}
\centering
\input{_L_9_1_2.tex}
\end{subfigure}
\hspace{5mm}
\begin{subfigure}[b]{0.3\textwidth}
\centering
\input{_L_9_1_3.tex}
\end{subfigure}
\caption{An upper bound for the $\L^{\star}$-invariant of $9_{1}$.}
\label{fig:L*_9_1}
\end{figure}

\clearpage

\begin{figure}[htbp]
\centering
\input{9_1_1-2.tex}
\caption{A tri-plane diagram of $9_{1}^{1,-2}$.}
\label{fig:9_1^1-2}
\end{figure}

\begin{figure}[htbp]
\centering
\begin{subfigure}[b]{0.3\textwidth}
\centering
\input{L_9_1_1-2_1.tex}
\end{subfigure}
\hspace{5mm}
\begin{subfigure}[b]{0.3\textwidth}
\centering
\input{L_9_1_1-2_2.tex}
\end{subfigure}
\hspace{5mm}
\begin{subfigure}[b]{0.3\textwidth}
\centering
\input{L_9_1_1-2_3.tex}
\end{subfigure}
\caption{An upper bound for the $\L$-invariant of $9_{1}^{1,-2}$.}
\label{fig:L_9_1^1-2}
\end{figure}

\begin{figure}[htbp]
\centering
\begin{subfigure}[b]{0.3\textwidth}
\centering
\input{_L_9_1_1-2_1.tex}
\end{subfigure}
\hspace{5mm}
\begin{subfigure}[b]{0.3\textwidth}
\centering
\input{_L_9_1_1-2_2.tex}
\end{subfigure}
\hspace{5mm}
\begin{subfigure}[b]{0.3\textwidth}
\centering
\input{_L_9_1_1-2_3.tex}
\end{subfigure}
\caption{An upper bound for the $\L^{\star}$-invariant of $9_{1}^{1,-2}$.}
\label{fig:L*_9_1^1-2}
\end{figure}

\clearpage

\begin{figure}[htbp]
\centering
\input{10_1_1.tex}
\caption{A tri-plane diagram of $10_{1}^{1}$.}
\label{fig:10_1^1}
\end{figure}

\begin{figure}[htbp]
\centering
\begin{subfigure}[b]{0.3\textwidth}
\centering
\input{L_10_1_1_1.tex}
\end{subfigure}
\hspace{5mm}
\begin{subfigure}[b]{0.3\textwidth}
\centering
\input{L_10_1_1_2.tex}
\end{subfigure}
\hspace{5mm}
\begin{subfigure}[b]{0.3\textwidth}
\centering
\input{L_10_1_1_3.tex}
\end{subfigure}
\caption{An upper bound for the $\L$-invariant of $10_{1}^{1}$.}
\label{fig:L_10_1^1}
\end{figure}

\begin{figure}[htbp]
\centering
\begin{subfigure}[b]{0.3\textwidth}
\centering
\input{_L_10_1_1_1.tex}
\end{subfigure}
\hspace{5mm}
\begin{subfigure}[b]{0.3\textwidth}
\centering
\input{_L_10_1_1_2.tex}
\end{subfigure}
\hspace{5mm}
\begin{subfigure}[b]{0.3\textwidth}
\centering
\input{_L_10_1_1_3.tex}
\end{subfigure}
\caption{An upper bound for the $\L^{\star}$-invariant of $10_{1}^{1}$.}
\label{fig:L*_10_1^1}
\end{figure}

\clearpage

\clearpage

\begin{figure}[htbp]
\centering
\input{10_1_001.tex}
\caption{A tri-plane diagram of $10_{1}^{0,0,1}$.}
\label{fig:10_1^001}
\end{figure}

\begin{figure}[htbp]
\centering
\begin{subfigure}[b]{0.3\textwidth}
\centering
\input{L_10_1_001_1.tex}
\end{subfigure}
\hspace{5mm}
\begin{subfigure}[b]{0.3\textwidth}
\centering
\input{L_10_1_001_2.tex}
\end{subfigure}
\hspace{5mm}
\begin{subfigure}[b]{0.3\textwidth}
\centering
\input{L_10_1_001_3.tex}
\end{subfigure}
\caption{An upper bound  for the $\L$-invariant of $10_{1}^{0,0,1}$.}
\label{fig:L_10_1^001}
\end{figure}

\begin{figure}[htbp]
\centering
\begin{subfigure}[b]{0.3\textwidth}
\centering
\input{_L_10_1_001_1.tex}
\end{subfigure}
\hspace{5mm}
\begin{subfigure}[b]{0.3\textwidth}
\centering
\input{_L_10_1_001_2.tex}
\end{subfigure}
\hspace{5mm}
\begin{subfigure}[b]{0.3\textwidth}
\centering
\input{_L_10_1_001_3.tex}
\end{subfigure}
\caption{An upper bound for the $\L^{\star}$-invariant of $10_{1}^{0,0,1}$.}
\label{fig:L*_10_1^001}
\end{figure}

\clearpage

\begin{figure}[htbp]
\centering
\input{10_1_-2-2.tex}
\caption{A tri-plane diagram of $10_{1}^{-2,-2}$.}
\label{fig:10_1^-2-2}
\end{figure}

\begin{figure}[htbp]
\centering
\begin{subfigure}[b]{0.3\textwidth}
\centering
\input{L_10_1_-2-2_1.tex}
\end{subfigure}
\hspace{5mm}
\begin{subfigure}[b]{0.3\textwidth}
\centering
\input{L_10_1_-2-2_2.tex}
\end{subfigure}
\hspace{5mm}
\begin{subfigure}[b]{0.3\textwidth}
\centering
\input{L_10_1_-2-2_3.tex}
\end{subfigure}
\caption{An upper bound for the $\L$-invariant of $10_{1}^{-2,-2}$.}
\label{fig:L_10_1^-2-2}
\end{figure}

\begin{figure}[htbp]
\centering
\begin{subfigure}[b]{0.3\textwidth}
\centering
\input{_L_10_1_-2-2_1.tex}
\end{subfigure}
\hspace{5mm}
\begin{subfigure}[b]{0.3\textwidth}
\centering
\input{_L_10_1_-2-2_2.tex}
\end{subfigure}
\hspace{5mm}
\begin{subfigure}[b]{0.3\textwidth}
\centering
\input{_L_10_1_-2-2_3.tex}
\end{subfigure}
\caption{An upper bound for the $\L^{\star}$-invariant of $10_{1}^{-2,-2}$.}
\label{fig:L*_10_1^-2-2}
\end{figure}

\end{document}